\documentclass{jT}
\usepackage{amssymb}
\usepackage{amsmath}
\usepackage{booktabs}
\usepackage{url}
\usepackage{pxfonts}
\usepackage{algorithmic}
\newcommand{\Z}{\mathbb{Z}}
\newcommand{\F}{\mathbb{F}}
\newcommand{\Hq}{\mathcal{H}_q}
\newcommand{\Fq}{\mathbb{F}_q}
\newcommand{\lcm}{{\rm lcm}}

\newtheorem{theorem}{Theorem}

\newtheorem{corollary}{Corollary}

\theoremstyle{definition}

\newtheorem{algorithm}{Algorithm}

\begin{document}

\title{On a theorem of Mestre and Schoof}
\author[John E. {\sc Cremona}]{{\sc John E.} CREMONA}
\address{John E. {\sc Cremona}\\
Mathematics Institute\\
University of Warwick\\
Coventry CV4 7AL\\
UK}
\email{J.E.Cremona@warwick.ac.uk}
\urladdr{\url{http://www.warwick.ac.uk/staff/J.E.Cremona/}}
\author[Andrew V. {\sc Sutherland}]{{\sc Andrew V. SUTHERLAND}}
\address{Andrew V. {\sc{Sutherland}}\\
Massachusetts Institute of Technology\\
Department of Mathematics\\
77 Massachusetts Avenue\\
Cambridge, MA 02139-4307\\
USA}
\email{drew@math.mit.edu}
\urladdr{\url{http://math.mit.edu/~drew/}}

\maketitle
\begin{resume}
Un th{\'e}or{\`e}me bien connu de Mestre et Schoof implique que la
cardinalit\'e d'une courbe elliptique $E$ d\'efinie sur un corps fini
$\Fq$ peut \^etre d\'etermin\'ee de mani\`ere univoque en calculant
les ordres de quelques points sur $E$ et sur sa tordue quadratique, \`a
condition que $q>229$. Nous \'etendons ce r\'esultat \`a tous les
corps finis avec $q>49$, et tous les corps premiers avec $q>29$.
\end{resume}
\begin{abstr}
A well known theorem of Mestre and Schoof implies that the order of an
elliptic curve $E$ over a prime field $\Fq$ can be uniquely determined
by computing the orders of a few points on $E$ and its quadratic twist,
provided that $q>229$.  We extend this result to all finite fields
with $q>49$, and all prime fields with $q>29$.
\end{abstr}
\bigskip
Let $E$ be an elliptic curve defined over the finite field $\Fq$ with~$q$
elements.  The number of points on $E/\Fq$, which we simply denote
$\#E$, is known to lie in the Hasse interval:
\begin{equation*}
\Hq = [q+1-2\sqrt{q},q+1+2\sqrt{q}].
\end{equation*}
Equivalently, the trace of Frobenius $t=q+1-\#E$ satisfies $|t|\le
2\sqrt{q}$.  A common strategy to compute $\#E$, when $q$ is not too
large, relies on the fact that the points on $E/\Fq$ form an abelian
group $E(\Fq)$ of order $\#E$.  For any $P\in E(\Fq)$, the integer
$\#E$ is a multiple of the order of $P$, and the multiples of~$|P|$
that lie in $\Hq$ can be efficiently determined using a baby-steps
giant-steps search.  If there is only one multiple in the interval,
it must be $\#E$; if not, we may try other $P\in E(\Fq)$ in the hope
of uniquely determining $\#E$.  This strategy will eventually succeed
if and only if the group exponent
$$\lambda(E)=\lcm\{|P|:P\in E(\Fq)\}$$ has a unique multiple in $\Hq$.
When this condition holds we expect to determine $\#E$ quite quickly:
with just two random points in $E(\Fq)$ we already succeed with
probability greater than $6/\pi^2$ (see \cite[Theorem
  8.1]{Sutherland:Thesis}).

Unfortunately, $\lambda(E)$ need not have a unique multiple in $\Hq$.
However, for prime~$q$ we have the following theorem of Mestre, as
extended by Schoof \cite[Theorem 3.2]{Schoof:ECPointCounting2}; the
result as stated in \cite{Schoof:ECPointCounting2} refers to the order
of a particular point $P$, but the following is an equivalent
statement.

\begin{theorem}[Mestre-Schoof]
Let $q>229$ be prime and $E$ an elliptic curve over $\Fq$ with
quadratic twist $E'$.  Either $\lambda(E)$ or $\lambda(E')$ has a
unique multiple in $\Hq$.
\end{theorem}

The quadratic twist $E'$ is an elliptic curve defined over~$\Fq$ that
is isomorphic to~$E$ over the quadratic extension~$\F_{q^2}$, and is
easily derived from $E$.  The orders of the groups~$E(\Fq)$
and~$E'(\Fq)$ satisfy~$\#E+\#E'=2(q+1)$.  For prime fields with
$q>229$, Theorem 1 implies that we may determine one of $\#E$ and
$\#E'$ by alternately computing the orders of points on $E$ and $E'$,
and once we know either $\#E$ or $\#E'$, we know both.

Theorem 1 does not hold for $q=229$.  Since there are counterexamples
whenever $q$ is a square, it does not hold in general for non-prime
finite fields either. The argument in the proof of~\cite[Theorem
3.2]{Schoof:ECPointCounting2} does not use the primality of~$q$, but
only that~$q$ is both large enough and not a square, so that the Hasse
bound on~$t$ cannot be attained.  If~$q=r^2$ is an even power of a
prime, then there are supersingular elliptic curves~$E$ over~$\Fq$
such that
$$
  E(\Fq)\cong \left(\Z/(r-1)\Z\right)^2\qquad\text{and}\qquad
  E'(\Fq)\cong \left(\Z/(r+1)\Z\right)^2.
$$ One may easily check that there are at least~$5$ multiples
of~$r-1$, and at least~$3$ multiples of~$r+1$, in~$\Hq$; however for
$r>7$ ($q>49$), the only pair that sum to~$2(q+1)$ are~$(r-1)^2$
and~$(r+1)^2$.  This resolves the ambiguity in these cases, leaving a
finite number of small exceptions.  For example, when $q=49$ there is
more than one pair of multiples of~$6$ and $8$ (respectively) which
sum to $2(q+1)=100$, since $100=36+64=60+40$.

The preceding observation led to this note, whose purpose is to extend
Theorem~1 to treat all finite fields (not just prime fields) $\Fq$
with $q>49$, and all prime fields with $q>29$.  Specifically, we prove
the following:

\begin{theorem}\label{theorem:CS}
Let $q\notin
\{3,4,5,7,9,11,16,17,23,25,29,49\}$ be a prime power, and let $E/\Fq$
be an elliptic curve.
Then there is a unique integer $t$ with $|t|\le 2\sqrt{q}$ such that
$\lambda(E)|(q+1-t)$ and $\lambda(E')|(q+1+t)$.
\end{theorem}

Our proof is entirely elementary, relying on just two properties of
elliptic curves over finite fields:

\begin{enumerate}
\item[(a)]
$\#E=q+1-t$ and $\#E'=q+1+t$ for some integer $t$ with $|t|\le 2\sqrt{q}$;
\vspace{2pt}
\item[(b)]
$E(\Fq)\cong \mathbb{Z}/n_1\mathbb{Z}\times \mathbb{Z}/n_2\mathbb{Z}$ with $n_1$ dividing both $n_2$ and $q-1$.
\end{enumerate}
Proofs of (a) and (b) may be found in most standard references,
including~\cite{Washington:EllipticCurves}.  We note that $n_2=\lambda(E)$,
and~$n_1=1$ when $E(\Fq)$ is cyclic.

\begin{proof}[Proof of Theorem 2]
Let $E$ be an elliptic curve over $\Fq$, and put $\#E=mM$ with
$M=\lambda(E)$, and $\#E'=nN$ with $N=\lambda(E')$. Without loss of
generality, we assume $a=q+1-\#E \ge 0$.  Taking $t=a$ shows existence, by
(a) and (b) above, so we need only prove that $t=a$ is the unique $t$
satisfying the conditions stated in the theorem.  For any such $t$ we
have $t\equiv q+1\bmod M$ and $t\equiv -(q+1)\bmod N$; hence $t$ lies
in an arithmetic sequence with difference $\lcm(M,N)$.  We also have
$|t|\le 2\sqrt{q}$; thus if $\lcm(M,N)>4\sqrt{q}$, then $t=a$ is
certainly unique.

We now show that $\lcm(M,N)\le 4\sqrt{q}$ implies $q\le 1024$.
We start from
\begin{equation*}
mMnN=(q+1-a)(q+1+a)=(q+1)^2-a^2\ge (q+1)^2-4q = (q-1)^2,
\end{equation*}
which yields
\begin{equation}\label{equation:one}
mn\ge\frac{(q-1)^2}{MN}=\frac{(q-1)^2}{\gcd(M,N)\lcm(M,N)}.
\end{equation}
Let $d=\gcd(m,n)$.  Then $d^2$ divides $mM+nN=2(q+1)$, so $d|(q+1)$, but also $d|(q-1)$, hence $d\le 2$. This implies $2\thinspace\lcm(M,N)\ge 2\thinspace\lcm(m,n)\ge mn$.  We also have $\gcd(M,N)\le\gcd(m,n)\gcd(M/m,N/n)\le 2\gcd(M/m,N/n)$.  Applying these inequalities to (\ref{equation:one}) we obtain
\begin{equation}\label{equation:two}
\lcm(M,N)^2\ge\frac{(q-1)^2}{4\gcd(M/m,N/n)}.
\end{equation}
We now suppose $\lcm(M,N)\le 4\sqrt{q}$, for otherwise the theorem holds.  We have $nN=q+1+a>q$, since we assumed $a\ge 0$, and $N\le 4\sqrt{q}$ implies that $n>\sqrt{q}/4$, so $N/n<16$.  Applying $\gcd(M/m,N/n)\le N/n<16$ to (\ref{equation:two}) yields
\begin{equation*}
4\sqrt{q} \ge \lcm(M,N)> (q-1)/8,
\end{equation*}
which implies that the prime power $q$ is at most $1024$.

The cases for $q\le 1024$ are addressed by a program listed in
the appendix that outputs the values of $q$, $M=\lambda(E)$, and
$N=\lambda(E')$ for which exceptions can arise.  This yields the set
of excluded $q$ and completes the proof.
\end{proof}

\textbf{Application.}  The proof of Theorem 2 suggests an
algorithm to compute~$\#E$, provided
that $q$ is small enough for the orders of randomly
chosen points in~$E(\Fq)$ to be easily computed.  It suffices to determine
integers $a$ and $m$ for which the set $S=\{x:x\equiv a
\bmod m\}$ contains $t=q+1-\#E$ but no $t'\ne t$ with $|t'|\le
2q$.  Beginning with $m=1$ and $a=0$, we compute $|P|$ for random points $P$
in $E(\Fq)$ or $E'(\Fq)$, and update $a$ and $m$ to reflect the fact
that $t\equiv q+1\bmod |P|$ when $P\in E(\Fq)$, and $t\equiv -(q+1)\bmod
|P|$ when $P\in E'(\Fq)$.  The new values of $a$ and $m$ may be
determined via the extended Euclidean algorithm.  When the set $S$
contains a unique~$t$ with $|t|\le 2\sqrt{q}$, we can conclude
that $\#E=q+1-t$ (and also that $\#E'=q+1+t$).

The probabilistic algorithm we have described is a \emph{Las Vegas}
algorithm, that is, its output is always correct and its expected
running time is finite.  The correctness of the algorithm follows from property (a). 
Theorem 2 ensures that the algorithm can terminate (provided that $q$ is not in the excluded set), 
and \cite[Theorem 8.2]{Sutherland:Thesis} bounds its expected running time.

An examination of Table 1 reveals that in many cases an ambiguous~$t'$ could be ruled out if $\lambda(E)$
or $\lambda(E')$ were known.  For example, when $q=49$, the trace $t'=-10$
yields $\#E=60$ and $\#E'=40$, so both $\lambda(E)$ and $\lambda(E')$ are
divisible by 5 (and are not 6 or 8).  If $E$ has trace $-10$, the
algorithm above will likely discover this and terminate within a few
iterations.  But when the trace of $E$ is 14 (and $\lambda(E)=6$ and $\lambda(E')=8$),
we can never be completely certain that we have ruled out $-10$ as a possibility.
Thus when an unconditional result is required, we must avoid
$q\in\{3,4,5,7,9,11,16,17,23,25,29,49\}$.

However, when $\lambda(E)$ and $\lambda(E')$ are known
we have the following corollary, which extends Proposition 4.19 of \cite{Washington:EllipticCurves}.

\begin{corollary}
Let $E/\Fq$ be an elliptic curve.  Up to isomorphism, the integers $\lambda(E)$
and $\lambda(E')$ uniquely determine the
groups~$E(\Fq)$ and $E'(\Fq)$, provided that $q\notin\{5,7,9,11,17,23,29\}$.  In every case,
$\lambda(E)$ and $\lambda(E')$ uniquely determine the set $\{E(\Fq),E'(\Fq)\}$.
\end{corollary}

Note that $\lambda(E)$ and $\#E$ together
determine $E(\Fq)$, by property (b).  To prove the corollary, apply Theorem 1 with 
a modified version of the algorithm in the appendix
that also requires $(q+1-t')/M$ to divide $M$ and $(q+1+t')/N$ to divide $N$.

As a final remark, we note that all the exceptional cases listed in
Table~\ref{table:Exceptions} can be eliminated if the orders of the
2-torsion and 3-torsion subgroups of $E(\Fq)$ are known (these orders
may be computed using the division polynomials).  
Alternatively, one can simply enumerate the points on $E/\Fq$ 
to determine $\#E$ when $q\le 49$.

\section*{Appendix}
For a prime power $q$, we wish to enumerate all $M$, $N$, and $t$ such that:
\begin{enumerate}
\item[(i)]
$M$ divides $q+1-t$ and $N$ divides $q+1+t$, with $0\le t\le 2\sqrt{q}$.
\item[(ii)]
$(q+1-t)/M$ divides $M$ and $q-1$, and $(q+1+t)/N$ divides $N$ and $q-1$.
\item[(iii)]
$M$ divides $q+1-t'$ and $N$ divides $q+1+t'$ for some $t'\ne t$ with $|t'|\le 2\sqrt{q}$.
\end{enumerate}
Any exception to Theorem \ref{theorem:CS} must arise from an elliptic
curve $E/\Fq$ with $\lambda(E)=M$, $\lambda(E')=N$, and $\#E=q+1-t$
(or from its twist, but the cases are symmetric, so we restrict to
$t\ge 0$).  Properties (i) and (ii) follow from (a) and (b) above, and
(iii) implies that $t$ does not uniquely satisfy the
requirements of the theorem.

Algorithm \ref{algorithm:enumerate} below finds all $M$, $N$, and $t$
satisfying (i), (ii), and (iii).  For $q\le 1024$, exceptional cases
are found only for the twelve values of~$q$ listed in Theorem~2.  Not
every case output by Algorithm 1 is actually realized by an elliptic
curve (in fact, all but one of the exceptions fail the condition that
$(q+1-t)/M\equiv(q+1+t)/N\pmod{2}$), but for each combination of $q$
and $t$ at least one is.  An example of each such case is listed in
Table \ref{table:Exceptions}, where we only list cases with $t\ge 0$:
for the symmetric cases with $t<0$, change the sign of $t$ and swap
$M$ and $N$.

\begin{algorithm}\label{algorithm:enumerate}
Given a prime power $q$, output all quadruples of integers $(M,N,t,t')$ satisfying (i), (ii), and (iii) above:
\end{algorithm}
\begin{algorithmic}
\FOR{all pairs of integers $(M,N)$ with $\sqrt{q}-1\le M,N\le 4\sqrt{q}$}
\FOR {all integers $t\in[0,2\sqrt{q}]$ with $M|(q+1-t)$ and $N|(q+1+t)$}
\STATE{Let $m=(q+1-t)/M$ and $n=(q+1+t)/N$}.
\IF{$m|M$ and $m|(q-1)$ and $n|N$ and $n|(q-1)$}
\FOR{all integers $t'\in[-2\sqrt{q},2\sqrt{q}]$}
\IF{$M|(q+1-t')$ and $N|(q+1+t')$}
\STATE{{\bf print} $M,N,t,t'$}.
\ENDIF
\ENDFOR
\ENDIF
\ENDFOR
\ENDFOR
\end{algorithmic}

\begin{table}[h]\label{table:exceptions}
\begin{center}
\begin{tabular}{@{}rrrrll@{}}
$q$ & $M$ & $N$ & $t$ & $E$&$t'$\\
\midrule
3&2&2&0&$y^2=x^3-x$&-2,2\\
4&1&3&4&$y^2+y=x^3+\alpha^2$&-2,1\\
5&2&4&2&$y^2=x^3+x$&-2\\
7&2&6&4&$y^2=x^3-1$&-2\\
7&4&4&0&$y^2=x^3+3x$&-4,4\\
9&2&4&6&$y^2 = x^3 + \alpha^2x$&-6,-2,2\\
11&4&8&4&$y^2=x^3+x+9$&-4\\
11&6&6&0&$y^2=x^3+2x$&-6,6\\
16&3&5&8&$y^2+y=x^3$&-7\\
17&6&12&6&$y^2=x^3 + x + 7$&-6\\
23&8&16&8&$y^2=x^3+5x+15$&-8\\
25&4&6&10&$y^2+y=x^3+\alpha^7$&-2\\
29&10&20&10&$y^2=x^3+x$&-10\\
49&6&8&14&$y^2=x^3+\alpha^2x$&-10\\
\bottomrule
\end{tabular}
\vspace{4pt}
\caption{Exceptional Cases with $t\ge 0$.}\label{table:Exceptions}
\small
The coefficient $\alpha$ denotes a primitive element of $\Fq$.
\end{center}
\end{table}

\bibliographystyle{amsplain}
\providecommand{\bysame}{\leavevmode\hbox to3em{\hrulefill}\thinspace}
\providecommand{\MR}{\relax\ifhmode\unskip\space\fi MR }
\providecommand{\MRhref}[2]{%
  \href{http://www.ams.org/mathscinet-getitem?mr=#1}{#2}
}
\providecommand{\href}[2]{#2}

\end{document}